\documentstyle[12pt,leqno]{amsart}
\newtheorem{theorem}{Theorem}[section]
\newtheorem{corollary}[theorem]{Corollary}

\newtheorem{remark}[theorem]{Remark}
\newtheorem{lemma}[theorem]{Lemma}

\begin{document}

\title[Prediction theory of stationary process]{A duality method in
   prediction theory of multivariate stationary sequences}
\author[M.~Frank and L.~Klotz]{{\rm By {\sc Michael Frank} and {\sc Lutz Klotz}
of Leipzig}}
\maketitle

%%%%%%%%%%%%%%%%%%%%%%%%%%%%%%%%%%%%%%%%%%%%%%%%%%%%%%%%%%%%%%%%%%%%%%%%%%%%%
\vspace{2cm}

\noindent
Kurztitel: Prediction theory of stationary sequences.

\vspace{1cm}


\vspace{1cm}


%%%%%%%%%%%%%%%%%%%%%%%%%%%%%%%%%%%%%%%%%%%%%%%%%%%%%%%%%%%%%%%%%%%%%%%%%%%%%

\vspace{2cm}

\noindent
{\bf Abstract.}
Let $W$ be an integrable positive Hermitian $q \times q$-matrix valued function
on the dual group of a discrete abelian group $\mathbf G$ such that $W^{-1}$
is integrable. Generalizing results of {\sc T.~Nakazi} and of {\sc
A.~G.~Miamee} \cite{N} and {\sc M.~Pourahmadi} \cite{MiP} for $q=1$ we
establish a correspondencebetween trigonometric approximation problems in
$L^2(W)$ and certain approximation problems in $L^2(W^{-1})$.
The result is applied to prediction problems for $q$-variate stationary
processes over $\mathbf G$, in particular, to the case ${\mathbf G} = \mathbb Z$.

\newpage

%\setcounter{equation}{0}

%%%%%%%%%%%%%%%%%%%%%%%%%%%%%%%%%%%%%%%%%%
% 0. Abschnitt:
\setcounter{equation}{0}
\renewcommand{\theequation}{\thesection.\arabic{equation}}

\section{{\bf Introduction}}

In 1984 {\sc T.~Nakazi} \cite{N} introduced a new idea into prediction theory
of univariate weakly stationary sequences. Under the additional assumption
that such a sequence has an absolutely continuous spectral measure and its
spectral density $w$ is such that $w^{-1}$ exists and is integrable
he related approximation problems in $L^2(w)$ to certain approximation
problems in $L^2(w^{-1})$. His method opened a way for him to give an elegant
proof of Szeg\"o's infimum formula and to obtain a partial solution of a
certain prediction problem which will be called Nakazi's prediction problem
in the present paper. {\sc A.~G.~Miamee} and {\sc M.~Pourahmadi} \cite{MiP}
pointed out that the essence of Nakazi's method consists in a certain duality
between the Hilbert spaces $L^2(w)$ and $L^2(w^{-1})$. This way they found a
unified approach to several `classical' prediction problems, and they
generalized some of Nakazi's results, partially even to more general
harmonizable stable sequences. The papers \cite{Mi,CMiP} contain further
completions of these results.

The aim of the present paper is the application of Nakazi's duality method to
$q$-variate ($q \in \mathbb N$, the set of positive integers) weakly
stationary processes over discrete abelian groups $\mathbf G$. Under the
assumption that the inverse $W^{-1}$ of the $q \times q$-matrix valued
spectral density $W$ of such a process is integrable we establish duality
relations between the left Hilbert modules $L^2(W)$ and $L^2(W^{-1})$.
Section three of our paper contains the general results. The further sections
are devoted to applications.

In Section four we obtain some well-known prediction results by quite different
proof methods. We compute the one-point interpolation error matrix and derive
Yaglom's interpolation recipe under the mentioned assumption on $W$. Moreover,
Section four contains our generalization of Nakazi's proof of Szeg\"o's infimum
formula to the multivariate situation. This way we obtain a special case of a
result due to {\sc V.~N.~Zasukhin} \cite{Z} as well as to {\sc H.~Helson} and
{\sc D.~Lowdenslager} \cite[Thm.~8]{HelL}.

If ${\mathbf G} = \mathbb Z$ is the abelian group of integers and the index
set $S$ of the known values consists of all negative
integers and the set of integers $\{ 1,2, ... ,n \}$ for some $n \in \mathbb
N$ then the arising prediction problem is called Nakazi's prediction problem.
In Section five we solve this problem under the additional assumption that
${\rm log} \, {\rm det} \, W$ is integrable.
We will see that in case $W^{-1}$ is integrable the
result is an easy consequence of our duality results. The more general case
in which merely ${\rm log} \, {\rm det} \, W$ is integrable can be solved by
approximation procedures. Section five also contains some straightforward
multivariate generalizations of univariate results of \cite{MiP}.

%%%%%%%%%%%%%%%%%%%%%%%%%%%%%%%%%%%%%
% 1. Abschnitt:
\setcounter{equation}{0}
\renewcommand{\theequation}{\thesection.\arabic{equation}}

\section{{\bf Preliminaries}}

Let $\mathbf G$ be a discrete abelian group with neutral element $0$,
${\mathbf G}^*$ be its dual group, and $\lambda$ be the normalized Haar
measure on ${\mathbf G}^*$, i.e.~$\lambda({{\mathbf G}^*})=1$. Relations
between measurable functions on ${\mathbf G}^*$ are to be understood as
relations which hold true almost everywhere (abbreviated to ``a.e.'') with
respect to (abbreviated to ``w.r.t.'') $\lambda$. Integration is always ment
to be done over ${\mathbf G}^*$. For a subset $S$ of $\mathbf G$ set $S^0 =
S \cup \{ 0 \}$ and $S^c = {\mathbf G} \setminus S^0$. For $q \in \mathbb N$
denote by ${\mathbf M}_q$ the algebra of all complex-valued $q \times q$-matrices.
The zero matrix and the identity matrix of ${\mathbf M}_q$ are denoted by
$0$ and $I$, respectively. For a non-empty subset $S$ of $\mathbf G$ let
${\mathcal T} (S)$ be the left ${\mathbf M}_q$-module of all
${\mathbf M}_q$-valued trigonometric polynomials with frequencies of $S$.
Functions of ${\mathcal T} (\{ 0 \})$, i.e.~constants, and elements of
${\mathbf M}_q$ will be identified and denoted by the same symbols.

If $A \in {\mathbf M}_q$ then the symbol $A^*$ stands for the adjoint of $A$,
${\rm det} \, A$ for its determinant, ${\rm tr} \, A$ for its normalized
trace and $|A|$ for its normalized euclidean norm, i.e.~$|A|^2 =
{\rm tr}(AA^*)$. If $A$ is regular then $A^{-1}$ denotes its inverse.
The set of all Hermitian matrices will be equipped with Loewner's
semi-ordering. In particular, a maximum or minimum of a subset of the set of
Hermitian $q \times q$-matrices is to be understood w.r.t.~that semi-ordering.
The cone of positive Hermitian $q \times q$-matrices is denoted by
${\mathbf M}_q^>$. For $A \in {\mathbf M}_q^>$ we denote its unique positive
square root by $A^{1/2}$. Finally, let $L^1$ and $L^2$ be the linear spaces
of (equivalence classes of) ${\mathbf M}_q$-valued functions on ${\mathbf
G}^*$ that are integrable or square-integrable, respectively. 

%%%%%%%%%%%%%%%%%%%%%%%%%%%%%%%%%%%%%%%%%%

Let ${\mathcal W}_q({{\mathbf G}^*})$ be the set of ${\mathbf M}_q^>$-valued
functions of $L^1$ and $\tilde{{\mathcal W}}_q({\mathbf G}^*) = \{ W \in
{\mathcal W}_q({{\mathbf G}^*}) : W^{-1} \in L^1 \}$.
For $W \in {\mathcal W}_q({{\mathbf G}^*})$ the symbol
$L^2(W)$ denotes the left Hilbert ${\mathbf M}_q$-module of (equivalence classes
of) ${\mathbf M}_q$-valued functions $F$ on ${\mathbf G}^*$ such that ${\rm
tr} (FWF^*)$ is integrable w.r.t.~$\lambda$. The ${\mathbf M}_q$-valued inner
product on $L^2(W)$ is defined by
\[
    (F,G) = \int FWG^* \, d\lambda \, ,
\]
the corresponding scalar product by
\[
    \langle F,G \rangle = {\rm tr} (F,G) \, , \: F,G \in L^2(W).
\]
The derived Hilbert space norm on $L^2(W)$ is denoted by $\| \cdot \|$. For a
more general construction see \cite{Ros}.

We recall that if $\mathcal M$ is a (closed) submodule of $L^2(W)$ and $F \in
L^2(W)$ there exists a unique function $F_{\mathcal M} \in \mathcal M$ such that
$(F-F_{\mathcal M}, F-F_{\mathcal M})$ is the minimum of the set $\{ (F-G,
F-G) : G \in {\mathcal M} \}$. The element $F_{\mathcal M}$ is the result of the
orthogonal projection of $F$ onto $\mathcal M$ w.r.t.~both $( \cdot , \cdot )$
and $\langle \cdot , \cdot \rangle$, i.e.~$(F-F_{\mathcal M}, G)=0$ for any
$G \in \mathcal M$. The orthogonal complement of $\mathcal M$ is equal to the
submodule $\{ H \in L^2(W) : (H,G)=0 \,\,\, {\rm for} \,\, {\rm all} \,\, G \in
\mathcal M \}$. The function $(F-F_{\mathcal M})$ is the result of the
orthogonal projection of $F$ onto the orthogonal complement of $\mathcal M$.
A detailed study of the geometry of Hilbert  ${\mathbf M}_q$-modules can be
found in \cite{GH}, for applications to prediction theory cf.~\cite{WMa}.

Considering $L^2(W)$ as the spectral domain of a $q$-variate weakly stationary
process over $\mathbf G$ we can formulate linear prediction problems of such
a process as trigonometric aopproximation problems in $L^2(W)$.
If $S \subseteq \mathbf G$ let ${\mathcal M} (S)$ be the closure of ${\mathcal T} (S)$
in $L^2(W)$ and $P_V$ be the orthogonal projection onto ${\mathcal M} (S)$. If $S
\subseteq {\mathbf G} \setminus \{ 0 \}$ the determination of $P_SI$ as well
as of the prediction error matrix
\begin{equation}     \label{2.1}
   \Delta_S = (I-P_SI, I-P_SI ) = (I-P_SI,I)
\end{equation}
is of fundamental importance in the prediction theory of $q$-variate weakly
stationary processes over $\mathbf G$. This kind of determination problem is
sometimes called the general prediction problem. There exist more or less
complete solutions of it for some special choices for the set $S$. For an
introduction to prediction theory of $q$-variate weakly stationary processes
we refer to \cite{Roz}.

%%%%%%%%%%%%%%%%%%%%%%%%%%%%%%%%%%%%%%%%%%%%%%%%%%%%%%%%%%%%%%%%%%%%%%%%%%%%%%%
% 2. Abschnitt:
\setcounter{equation}{0}

\section{{\bf Duality relations}}

Let $W \in \tilde{{\mathcal W}_q}({{\mathbf G}^*})$. Then along with $L^2(W)$
one can consider the left Hilbert ${\mathbf M}_q$-module $L^2(W^{-1})$.
The symbols $(\cdot ,\cdot )$,
$\langle \cdot ,\cdot \rangle$, $\| \cdot \|$, ${\mathcal M}(S)$, $P_S$ and
$\Delta_S$ referring to $L^2(W)$ will be replaced by $(\cdot ,\cdot )_\sim$,
$\langle \cdot ,\cdot \rangle_\sim$, $\| \cdot \|_\sim$,
$\tilde{{\mathcal M}}(S)$, $\tilde{P}_S$ and $\tilde{\Delta}_S$, respectively, for
the denotation of the corresponding objects related to $L^2(W^{-1})$. It turns
out that a certain duality between $L^2(W)$ and $L^2(W^{-1})$ is helpful if
we study the general prediction problem in $L^2(W)$.

\begin{lemma}   \label{Lemma3.1}
    Let $W \in \tilde{{\mathcal W}_q}({{\mathbf G}^*})$. Then the mapping
    \[
       F \longrightarrow FW \,\, , \,\, F \in L^2(W) \, ,
    \]
    is an isometric isomorphism of $L^2(W)$ onto $L^2(W^{-1})$.
\end{lemma}

P r o o f :
If $F,G \in L^2(W)$, we obviously have $(F,G) = (FW,GW)_\sim$. If $H
\in L^2(W^{-1})$ then $HW^{-1} \in L^2(W)$ and $HW^{-1}W = H$.
\hfill $\Box$

Another immediate consequence of the integrability of $W^{-1}$ is the point of
view on $L^2(W)$ and on $L^2(W^{-1})$ as on subspaces of $L^1$. This can be seen
from the inequality
\begin{eqnarray}     \label{3.1}
   \int |F| \, d \lambda = \int |FW^{1/2}W^{-1/2}| \, d \lambda
       & \leq &
       \int | FW^{1/2}| \cdot | W^{-1/2}| \, d \lambda \\
       & \leq &
       \left( \int |FW^{1/2}|^2 \, d \lambda \right)^{1/2}
       \left( \int | W^{-1}| \, d \lambda \right)^{1/2}  \nonumber\\
       & < &
       \infty \,\, , \,\, F \in L^2(W) \, .   \nonumber
\end{eqnarray}
If $\mathcal M$ is a submodule of $L^2(W)$ then the submodule
\[
    {\mathcal M}^d = \left\{ G \in L^2(W^{-1}) \, : \, \int FG^* \, d \lambda
    = 0 \, , \, F \in {\mathcal M} \right\} \subseteq L^2(W^{-1})
\]
is said to be its dual. Note, that ${\mathcal M}^d$ is the orthogonal complement
in $L^2(W^{-1})$ of the submodule ${\mathcal M}W = \{ FW : F \in \mathcal M \}$.
For a subset $S \subseteq {\mathbf G} \setminus \{ 0 \}$ the symbol $\tilde{I}_S$
denotes the orthogonal projection of $I$ onto ${\mathcal M}(S^0)^d \subseteq
L^2(W^{-1})$.

\begin{theorem}   \label{Theorem3.2}
     Let $W \in \tilde{{\mathcal W}_q}({{\mathbf G}^*})$. Then for subsets $S
     \subseteq {\mathbf G} \setminus \{ 0 \}$
     \begin{equation}  \label{3.2}
        P_SI = I - (I -\tilde{I}_S, I)_\sim^{-1} (I-\tilde{I}_S)W^{-1}
     \end{equation}
     and
     \begin{equation}     \label{3.3}
        \Delta_S = (I-\tilde{I}_S,I)_\sim^{-1} = (I-\tilde{I}_S,I-\tilde{I}_S
        )_\sim^{-1} \, .
     \end{equation}
\end{theorem}

In order to prove Theorem \ref{Theorem3.2} we need the following lemma.

\begin{lemma} \label{Lemma3.3}
     Let $F \in {\mathcal M}(S^0)$. Then $F \in {\mathcal M} (S)$ if and only if
     $\int F \, d \lambda = 0$.
\end{lemma}

P r o o f :
If $F \in {\mathcal T} (S)$ then $\int F \, d \lambda = 0$. If $F \in {\mathcal M} (S)$
we approximate it by elements of ${\mathcal T} (S)$, and by (\ref{3.1}) we obtain
$\int F \, d \lambda = 0$. If $F \in {\mathcal M}(S^0)$ it can be written in the
form $F=A+G$, where $A \in \mathcal M$ and $G \in {\mathcal M} (S)$. By the result
just proved we have $\int F \, d \lambda =A$. Consequently, $\int F \, d \lambda
= 0$ implies $F=G \in {\mathcal M} (S)$.
\hfill $\Box$

\medskip
P r o o f  $\,$  o f $\,$ T h e o r e m $\,$ \ref{Theorem3.2} :
Since $\tilde{I}_S \in {\mathcal M}(S^0)^d$ we have $\int \tilde{I}_S \, d
\lambda = \int \tilde{I}_S I^* \, d \lambda = 0$ and hence, $(I-\tilde{I}_S,W)_\sim
= \int (I-\tilde{I}_S) \, d \lambda =I$. By a generalization of the Cauchy-Schwarz
inequality the inequality $(W,W)_\sim^{-1} = (I-\tilde{I}_S,W)_\sim
(W,W)_\sim^{-1} (I-\tilde{I}_S,W)_\sim \leq (I-\tilde{I}_S, I-\tilde{I}_S)_\sim$
holds, cf.~\cite[Thm.~1]{B}. Therefore the right hand side of (\ref{3.2}) is
well-defined. Let us denote it by $F_S$. Since $(I-\tilde{I}_S) \in {\mathcal
M}(S^0)W$ we conclude $F_S \in {\mathcal M}(S^0)$ by Lemma \ref{Lemma3.1}.
Because $\int F_S \, d \lambda = I-(I-\tilde{I}_S,I)_\sim^{-1} (I-\tilde{I}_S,I)_\sim
=0$ Lemma \ref{Lemma3.3} implies
\begin{equation} \label{3.4}
   F_S \in {\mathcal M}(S) \, .
\end{equation}
Let $F \in {\mathcal M}(S)$. Then $(I,FW)_\sim = \int F^* \, d \lambda = 0$ by
Lemma \ref{Lemma3.3}. Therefore, the equality $(I-F_S,F) = (I-\tilde{I}_S,I)_\sim^{-1}
(I-\tilde{I}_S,FW)_\sim = (I-\tilde{I}_S,I)_\sim^{-1} (I,FW)_\sim = 0$ yields.
By a combination with (\ref{3.4}) we get $F_S=P_SI$. To obtain formula (\ref{3.3})
simply replace one term according to (\ref{3.2}) in (\ref{2.1}).
\hfill $\Box$

\medskip
From Theorem \ref{Theorem3.2} we can derive a general result which seems not
have immediate applications to prediction theory of multivariate stationary
processes, however which can be considered as a certain generalization of the
univariate case.

Let ${\mathbf M}_q'$ be the set of all $q \times q$-matrices whose determinant
equals to one. For $S \subseteq {\mathbf G} \setminus \{0\}$ define
\begin{equation} \label{3.5}
    \delta_S = \inf \{ \| A-F \|^2 : A \in {\mathbf M}_q' \, \, , \, \, F \in
    {\mathcal M}(S) \} \, .
\end{equation}

\begin{lemma}  \label{Lemma3.4}
   Let $W \in {\mathcal W}_q({{\mathbf G}^*})$ and $S$ be a subset of ${\mathbf G}
   \setminus \{ 0 \}$. Then
   \begin{equation} \label{3.6}
       \delta_S = [ {\rm det}(\Delta_S) ]^{1/q} \, .
   \end{equation}
\end{lemma}

P r o o f : We have the equality $[{\rm det} (I-F,I-F)]^{1/q} = \inf \{ {\rm tr}
((I-F,I-F)B) : B \in {\mathbf M}_q^> \cap {\mathbf M}_q') \}$, $F \in {\mathcal M}(S)$,
cf.~\cite[Section 7.8, Problem 19]{HoJ}. Since $B=A^*A$ for some $A \in {\mathcal M}'$ we
obtain $[{\rm det} (I-F,I-F)]^{1/q} = \inf \{ \| A-AF \|^2 : A \in {\mathcal M}' \}$,
$F \in {\mathcal M}(S)$, and therefore
\begin{eqnarray*}
   [{\rm det} (\Delta_S)]^{1/q} & = & \{ {\rm det} [ \min \{ (I-F,I-F) : F \in
   {\mathcal M}(S) \} ] \}^{1/q} \\
   & = & \min \{ [ {\rm det} (I-F,I-F)]^{1/q} : F \in {\mathcal M}(S) \} \\
   & = & \inf \{ \| A-AF \|^2 : A \in {\mathbf M}_q' \,\, , \,\, F \in
   {\mathcal M}(S) \} \\
   & = & \delta_S  \, .
\end{eqnarray*}
$\,$ \hfill $\Box$

\begin{remark} \label{Remark3.5}
   {\rm
   Since a matrix $A \in {\mathbf M}_q'$ can be written in the form $A=UB$
   for a unitary matrix $U$ and for a matrix $B \in ({\mathbf M}_q^> \cap
   {\mathbf M}_q')$ the value of $\delta_S$ does not change if $A$ runs through
   $({\mathbf M}_q^> \cap {\mathbf M}_q')$ on the right-hand side of (\ref{3.5}).
   Moreover, if $A$ exhausts the larger set of all $q \times q$-matrices with
   $| {\rm det} \, A|=1$ we get the same result.
   }
\end{remark}

Combining Theorem \ref{Theorem3.2} and Lemma \ref{Lemma3.4} we get the following
assertion.

\begin{theorem} \label{Theorem3.6}
    Let $W \in\tilde{{\mathcal W}}_q({{\mathbf G}^*})$. Then for subsets
    $S \subseteq {\mathbf G} \setminus \{ 0 \}$ we have the equality
    \begin{equation} \label{3.7}
       \delta_S = ( \inf \{ \| A-G \|_\sim^2 : A \in {\mathbf M}_q' \, \, , \, \,
       G \in {\mathcal M}(S^0)^d \} )^{-1} \, .
    \end{equation}
\end{theorem}

P r o o f :
In a similar way to that one taken in the proof of Lemma \ref{Lemma3.4} we
derive the equality
\begin{eqnarray*}
   \inf \{ \| A-G \|_\sim^2 : A \in {\mathbf M}_q' \, \, , \, \, G \in {\mathcal M}
   (S^0)^d \} & = &
   \{ {\rm det} [ \min \{ (I-G,I-G)_\sim : G \in {\mathcal M}(S^0)^d \} ] \}^{1/q}
   \\
   & = & [ {\rm det} (I-\tilde{I}_S,I-\tilde{I}_S)_\sim]^{1/q} \, .
\end{eqnarray*}
By (\ref{3.3}) and (\ref{3.6}) the right-hand end of this equality equals to
$\delta_S^{-1}$.
\hfill $\Box$

Another duality relation can be derived by a simple adaptation of a chain of
equalities mentioned in the proof of \cite[Thm.~4.1]{C} to the
multivariate situation.

\begin{theorem} \label{Theorem3.7}
    Let $W \in \tilde{{\mathcal W}}_q ({{\mathbf G}^*})$. Then for subsets $S
    \subseteq {\mathbf G}\setminus \{ 0 \}$ one has
    \begin{equation}  \label{3.8}
         \inf \{ \| A-P_SI \| : A \in {\mathbf M}_q \,\, , \,\, {\rm tr} \, A =1
         \} = \| I-\tilde{I}_S \|_\sim^{-1} \, .
    \end{equation}
\end{theorem}

P r o o f :
We have the following chain of equalities
\begin{eqnarray*}
    \lefteqn{
    \inf \{ \| A-P_SI \|^2 : A \in {\mathbf M}_q \,\, , \,\, {\rm tr} \, A =1 \} =
    } \\
    & = &
    \inf \{ \| A-T \|^2 : A \in {\mathbf M}_q \,\, , \,\, {\rm tr} \, A =1 \,\, ,
    \,\, T \in {\mathcal T}(S)\}  \\
    & = &
    \inf \left\{ \frac{\| A-T \|^2}{| {\rm tr} \, A |^2} : A \in {\mathbf M}_q \,\, ,
    \,\, {\rm tr} \, A \not= 0 \,\, , \,\, T \in {\mathcal T}(S) \right\} \\
    & = &
    \inf \left\{ \frac{\| T \|^2}{| {\rm tr} ( \int T \, d \lambda ) |^2} :
    T \in {\mathcal T}(S^0) \, \, , \, \, {\rm tr} \left(\int T \, d \lambda
    \right) \not= 0 \right\} \\
    & = &
    \left( \sup \left\{ \frac{\langle TW,I \rangle_\sim}{\| TW \|^2_\sim} :
    T \in {\mathcal T}(S^0) \, \, , \, \, {\rm tr} \left(\int T \, d \lambda
    \right) \not= 0 \right\} \right)^{-1}\\
    & = &
    \| \tilde{P}_{S^0} I \|_\sim^{-2} \\
    & = &
    \| I-\tilde{I}_S \|_\sim^{-2} \, . 
\end{eqnarray*}
$\,$ \hfill $\Box$

\medskip
For applications of the preceding theorems a description of the space
${\mathcal M}(S^0)^d$ is needed. The identification of ${\mathcal M}(S^0)^d$
and of $\tilde{{\mathcal M}}(S^c)$ would be desirable from the point of view
of prediction theory. Clearly,
\begin{equation} \label{3.9}
    \tilde{{\mathcal M}}(S^c) \subseteq {\mathcal M}(S^0)^d  \, .
\end{equation}
However, whether equality really holds or not,
seems to be a difficult problem, in general. It is related to basis properties
of characters. {\sc Miamee} and {\sc Pourahmadi} \cite{MiP} discussed this
question and gave particular answers to it in case $q=1$, cf.~also \cite{Mi}.
In contrast to these careful investigations, in the proofs of \cite[Thm.~4.1]{C}
and of \cite[Thm.~1]{CMiP} the equality $\tilde{{\mathcal M}}(S^c) =
{\mathcal M}(S^0)^d$ seems to be used in a rather general situation, but
without any explanation how to prove it.

We call a subset $S \subseteq {\mathbf G}$ to be ${\mathbf G}$-exact if for any $q \in
\mathbb N$ and any $W \in \tilde{{\mathcal W}}_q({{\mathbf G}^*})$ the set
identity $\tilde{{\mathcal M}}({\mathbf G} \setminus S) = {\mathcal M}(S)^d$
holds. Because of the symmetry of this definition we immediately obtain the
following fact.

\begin{lemma} \label{Lemma3.8}
   The set $S$ is $\mathbf G$-exact if and only if the set ${\mathbf G} \setminus S$
   is $\mathbf G$-exact.
\end{lemma}

The next result demonstrates that the alteration of a $\mathbf G$-exact set
$S$ by finitely many elements does not have any influence on its
$\mathbf G$-exactness.

\begin{theorem} \label{Theorem3.9}
   Let $S$ be a subset of $\mathbf G$ and $g \in {\mathbf G} \setminus S$.
   Then $S$ is $\mathbf G$-exact if and only if $S \cup \{ g \}$ is $\mathbf
   G$-exact.
\end{theorem}

P r o o f :
Without loss of generality we may assume that $0 \not\in S$ and $g=0$.
Suppose $S$ to be $\mathbf G$-exact. By (\ref{3.9}) we have the chain of set
inclusions $\tilde{{\mathcal M}}(S^c) \subseteq {\mathcal M}(S^0)^d \subseteq
{\mathcal M}(S)^d \subseteq \tilde{{\mathcal M}}({\mathbf G} \setminus S)$.
Thus, if $G \in {\mathcal M}(S^0)^d$ it can be decomposed as $G=A+H$ for
some $A \in {\mathbf M}_q$ and $H \in \tilde{{\mathcal M}}(S^c)$. So $G-H =A
\in {\mathcal M}(S^0)^d$ and hence, $A=0$ and $G=H \in \tilde{{\mathcal M}}(S^c)$.
The latter means that $S^0$ is $\mathbf G$-exact. To derive the $\mathbf G$-exactness
of $S$ from the $\mathbf G$-exactness of $S^0$ one has to combine the result
just proved and Lemma \ref{Lemma3.8}.
\hfill $\Box$

\begin{corollary}  \label{Corollary3.10}
   If either $S$ or ${\mathbf G} \setminus S$ are finite subsets of $\mathbf G$
   then $S$ is $\mathbf G$-exact.
\end{corollary}

P r o o f :
Define ${\mathcal M}(\emptyset) = \{ 0 \}$, so the empty set is $\mathbf G$-exact
by definition. With this setting the assertion follows from Theorem \ref{Theorem3.9}
and Lemma \ref{Lemma3.8}.
\hfill $\Box$

\medskip
If ${\mathbf G} = \mathbb Z$ then its dual group ${\mathbf G}^*$ can be
identified with the numerical interval $[ -\pi,\pi )$, where the addition of
elements is understood to be done ${\rm mod } \, 2\pi$. The characters of
${{\mathbf G}^*} = [-\pi,\pi)$ can be described as the set of functions
$e_k(t) = {\rm e}^{ikt}$, $t \in [-\pi,\pi)$, $k \in \mathbb Z$. A univariate
version of the assertion stated thereafter was given in \cite[Lemma 3.6]{MiP},
compare with \cite[Thm.~3.1]{Mi} for a univariate extension of this
assertion to $L^p$-spaces.

\begin{theorem}  \label{Theorem3.11}
     The set $\mathbb N$ is $\mathbb Z$-exact.
\end{theorem}

P r o o f :
By (\ref{3.9}) we have $\tilde{{\mathcal M}}({\mathbb Z} \setminus {\mathbb N})
\subseteq {\mathcal M}({\mathbb N})^d$. Let $G \in {\mathcal M}({\mathbb N})^d$
be orthogonal to $\tilde{{\mathcal M}}({\mathbb Z} \setminus {\mathbb N})$.
Then
\begin{equation} \label{3.10}
     \int e_k G^* \, d \lambda = 0 \,\, , \,\, k \in \mathbb N \, ,
\end{equation}
and
\begin{equation}  \label{3.11}
     \int e_k W^{-1}G^* \, d \lambda = 0 \, \, , \, \, k \in {\mathbb Z}
     \setminus \mathbb N \, .
\end{equation}
By (\ref{3.10}) the element $G^*$ belongs to the Hardy space ${\mathbf H}^1$
(of ${\mathbf M}_q$-valued functions), and (\ref{3.11}) implies that $GW^{-1}$
belongs to ${\mathbf H}^1$, too. Therefore, $GW^{-1}G^*$ is an element of
${\mathbf H}^{1/2}$ taking values in the set of non-negative Hermitian matrices.
According to \cite{ST} the function $GW^{-1}G^*$ is a constant function.
However, by (\ref{3.10}) and (\ref{3.11}) the index zero Fourier coefficient
of $GW^{-1}G^*$ equals to zero, so $GW^{-1}G^* = 0$ and $G=0$.
\hfill $\Box$

%%%%%%%%%%%%%%%%%%%%%%%%%%%%%%%%%%%%%%%%%%%%%%%%%%%%%%%%%%%%%%%%%%%%%%%%%%%%%%%
% 3. Abschnitt:
\setcounter{equation}{0}

\section{{\bf Some applications}}

From the general assertions of section three we can easily derive some prediction
results for multivariate stationary sequences. Many of them are well-known, however
here they are drawn from a different context.

\begin{theorem} {\rm (cf.~\cite[Thm.~4.3]{Heb} or \cite[Cor.~2.9]{Ma})} \label{Theorem4.1}
     \newline
     Let $W \in \tilde{{\mathcal W}}_q({{\mathbf G}^*})$ and $S= {\mathbf G}
     \setminus \{ 0 \}$. Then
     \begin{eqnarray*}
     P_SI & = & I - \left( \int W^{-1} \, d \lambda \right)^{-1} W^{-1}  \, , \\
     \Delta_S & = & \left( \int W^{-1} \, d \lambda \right)^{-1}   \,\,\, {\rm and} \\
     \delta_S & = & \left[ {\rm det} \left( \int W^{-1} \, d \lambda \right)\right]^{-1/q} \, .
     \end{eqnarray*}
\end{theorem}

P r o o f :
Since ${\mathcal M}(S^0)^d= \{ 0 \}$ the first two results are immediate
consequences of Theorem \ref{Theorem3.2}. The third result can be derived from
(\ref{3.6}). \hfill $\Box$

\medskip
Now let ${\mathbf G}= \mathbb Z$, and for a certain $n \in \mathbb N$ set
\[
    S_1 = {\mathbb Z} \setminus \{ 0,1,2, ...,n\} \, .
\]

\begin{theorem}  {\rm (cf.~\cite[Thm.~2]{Ya} for the univariate case)} \label{Theorem4.2}
     \newline
     Let $n \in \mathbb N$, $W \in \tilde{{\mathcal W}}_q([-\pi,\pi))$ and
     $C_l = \int e_lW^{-1} \, d \lambda$, $l \in \mathbb Z$. Let the $q \times
     q$-matrices $\{ D_k : k=1,...,n \}$ be the solutions of the linear system
     \begin{equation}  \label{4.1}
         \sum_{k=1}^n C_{j-k}D_k^* = C_j \,\, , \,\, j=1,...,n \, .
     \end{equation}
     Then
     \begin{eqnarray}   \label{4.2}
          P_{S_1}I & = & I - \left( I - \sum_{k=1}^n D_ke_k,I \right)_\sim^{-1}
          \left[ I-\sum_{k=1}^n D_ke_k) W^{-1} \right] \, , \nonumber\\
          \Delta_{S_1} & = & \left( I - \sum_{k=1}^n D_ke_k,I \right)_\sim^{-1}
     \end{eqnarray}
     and
     \begin{eqnarray*}
          \delta_{S_1} & = & \left[ {\rm det} \left( I - \sum_{k=1}^n D_ke_k,I
          \right)_\sim  \right]^{-1/q} \, . 
     \end{eqnarray*}
\end{theorem}

P r o o f :
According to Corollary \ref{Corollary3.10} the set $S_1$ is $\mathbb Z$-exact.
Elementary calculations show that $\tilde{P}_{S_1^c} I = \sum_{k=1}^n D_ke_k$.
Knowing this the result is an immediate consequence of Theorem \ref{Theorem3.2}
and Lemma \ref{Lemma3.4}.
\hfill $\Box$

\begin{remark}   {\rm
From (\ref{4.2}) we obtain
\begin{equation}  \label{4.3}
    \sum_{k=1}^n C_{-k}D_k^* = C_0 - \Delta_{S_1}^{-1}  \, .
\end{equation}
Let $\Gamma_n$ be the Toeplitz matrix $\Gamma_n = ( C_{j-k} )_{j,k=1,...,n}$,
$n \in \mathbb N$. It is not hard to see that the system of equations (\ref{4.1})
and (\ref{4.3}) can be written in the form
\[
    \Gamma_{n+1} \cdot \left( \begin{array}{c} -I \\ D_1^* \\ . \\ . \\ . \\ D_n^*
    \end{array} \right) =
    \left( \begin{array}{c} -\Delta_{S_1}^{-1} \\ 0 \\ . \\ . \\ . \\ 0
    \end{array} \right)  \,\, , \quad {\rm i.e.} \quad
    \left( \begin{array}{c} -I \\ D_1^* \\ . \\ . \\ . \\ D_n^*
    \end{array} \right) =  \Gamma_{n+1}^{-1} \cdot
    \left( \begin{array}{c} -\Delta_{S_1}^{-1} \\ 0 \\ . \\ . \\ . \\ 0
    \end{array} \right)  \,\, .
\]
Comparing the first $q \times q$-blocks from both the sides we get $I =
\Omega_n \Delta_{S_1}^{-1}$, and therefore $1= {\rm det} ( \Omega_n) {\rm det}
(\Delta_{S_1}^{-1})$. Here $\Omega_n$ denotes the $q \times q$-matrix in the
left upper corner of the matrix $\Gamma_{n+1}^{-1}$. Using a well-known rule
for computing minors of inverse matrices together with the fact that the
$nq \times nq$-matrix in the right lower corner of $\Gamma_{n+1}$ equals
$\Gamma_n$, we get
\[
    {\rm det} (\Delta_{S_1} ) =
    {\rm det} (\Omega_n)  =
    \frac{{\rm det} (\Gamma_n)}{{\rm det} (\Gamma_{n+1})}
\]
for the determinant of the prediction error matrix.

\medskip
For the special case
where $q=1$ {\sc Nakazi} gave an elegant proof of Szeg\"o's infimum formula,
cf.~\cite[Cor.~4]{N}. {\sc Zasukhin} \cite{Z}, {\sc Helson} and {\sc Lowdenslager}
\cite[Thm.~8]{HelL} stated a version of Szeg\"o's theorem for non-negative
Hermitian matrix-valued measures. Using the results of section three we can adapt Nakazi's
proof to the multivariate situation to obtain the result by Zasukhin and Helson-Lowdenslager
for absolutely continuous non-negative Hermitian matrix-valued measures.
}
\end{remark}

\begin{theorem} \label{Theorem4.3}
    Let ${\mathbf G} = \mathbb Z$ and $W$ be a function of $L^1$, whose
    values are non-negative Hermitian matrices. Then
    \begin{equation}      \label{4.4}
         \delta_{\mathbb N} = {\rm exp} \left\{ \frac{1}{q} \cdot \int
         {\rm log} \, {\rm det} \, W \, d \lambda \right\} \, ,
    \end{equation}
    where the right-hand side of (\ref{4.4}) has to be interpreted as zero
    in case the expression ${\rm log} \, {\rm det} \, W$ is not integrable.
\end{theorem}

P r o o f :
First assume that $W \in \tilde{{\mathcal W}}_q([-\pi,\pi))$. If $A \in
{\mathbf M}_q'$ and $T \in {\mathcal T}({\mathbb N})$ from the inequality
between the arithmetic and the geometric means and from the fact that
$\int {\rm log} \, |{\rm det} (A-T)| \, d \lambda \geq {\rm log} \,
|{\rm det} \, A| = 0$ we obtain
\begin{eqnarray*}
     \| A-T \|^2 & = & \int {\rm tr} [(A-T)W(A-T)^*] \, d \lambda \\
     & \geq &
     {\rm exp} \left\{ \int {\rm log} \, {\rm tr} [(A-T)W(A-T)^*] \, d \lambda
     \right\}\\
     & \geq & 
     {\rm exp} \left\{ \int {\rm log} \left\{ {\rm det} [(A-T)W(A-T)^*]
     \right\}^{1/q} \, d \lambda \right\}\\
     & \geq &
     {\rm exp} \left\{ \frac{1}{q} \cdot \int {\rm log} \, {\rm det}  \, W \, d \lambda
     \right\} \, .
\end{eqnarray*}
Hence
\begin{equation} \label{4.5}
   \delta_{\mathbb N} \geq
   {\rm exp} \left\{ \frac{1}{q} \cdot \int {\rm log} \, {\rm det} \, W \, d \lambda
     \right\} \, .
\end{equation}
Similarly,
\begin{equation} \label{4.6}
   \tilde{\delta}_{{\mathbb N}^c} \geq
   {\rm exp} \left\{ \frac{1}{q} \cdot \int {\rm log} \, {\rm det} (W^{-1}) \, d \lambda
     \right\} \, .
\end{equation}
By Theorem \ref{Theorem3.11} and Theorem \ref{Theorem3.9} the set ${\mathbb N}^0$ is
$\mathbb Z$-exact, and by (\ref{3.7}) the equality $\delta_{\mathbb N} =
\tilde{\delta}_{{\mathbb N}^c}^{-1}$ turns out. Comparing this with (\ref{4.5}) and
(\ref{4.6}) we get the desired result under the additional assumption that
$W^{-1} \in L^1$. The general assertion can be derived from this partial one
by a standard approximation argument demonstrated in \cite[pp.~260-261]{A}
for the case $q=1$.
\hfill $\Box$

\medskip
From Theorem \ref{Theorem4.3} we obtain a somewhat surprising result.

\begin{corollary}  \label{Corollary4.4}
  Let ${\mathbf G} = \mathbb Z$ and $W \in {\mathcal W}_q([-\pi,\pi))$.
  Then
  \[
   \left[ {\rm det} (\Delta_{\mathbb N} ) \right]^{1/q} =
   \inf \left\{ \int [ {\rm det} [(I-T)W(I-T)^*]]^{1/q} \, d\lambda \, : \,
   T \in {\mathcal T}({\mathbb N}) \right\} \, .
  \]
\end{corollary}

P r o o f :
By (\ref{3.6}) and (\ref{4.4}) we conclude $\left[ {\rm det}
(\Delta_{\mathbb N} ) \right]^{1/q} =  {\rm exp} \left\{ 1/q \cdot \int {\rm log}
\, {\rm det} \, W \, d \lambda \right\}$.  On the other hand,
\[
 {\rm exp} \left\{ 1/q \cdot \int {\rm log} \, {\rm det} \, W \, d \lambda
 \right\} = \inf \left\{ \int [ {\rm det} [(I-T)W(I-T)^*]]^{1/q} \, d\lambda
 \, : \, T \in {\mathcal T}({\mathbb N}) \right\} \, ,
\]
cf.~\cite[Cor.]{K}.
\hfill $\Box$

%%%%%%%%%%%%%%%%%%%%%%%%%%%%%%%%%%%%%%%%%%%%%%%%%%%%%%%%%%%%%%%%%%%%%%%%%%%%%%%
% 4. Abschnitt:
\setcounter{equation}{0}

\section{{\bf Nakazi's prediction problem}}

Another application of the results of Section three gives some information about the
multivariate version of Nakazi's prediction problem. Let ${\mathbf G} = \mathbb Z$
and $S$ be of the form
\begin{equation} \label{2.2}
   S_2 = {\mathbb N}^c \cup \{ 1,2,...,n \}
\end{equation}
for some $n \in \mathbb N$.
Assume the function $W$ belongs to ${\mathcal W}_q([-\pi,\pi ))$ and is such
that ${\rm log} \, {\rm det} \, W$ is integrable. The set of such functions
is denoted by ${\mathcal W}_q'([-\pi,\pi])$ in the sequel. Let $\Phi$ be the
unique outer function of the Hardy space ${\mathbf H}^2$ (of functions with
values in ${\mathbf M}_q$) such that
\begin{equation}  \label{5.1}
   W= \Phi^*\Phi \qquad {\rm and} \qquad \int \Phi \, d \lambda \in
   {\mathbf M}_q^> \, .
\end{equation}
The function $\Phi^{-1}$ is also outer, but it does not belong to
${\mathbf H}^2$, in general. Let $B_j$ be the $j$-th Taylor coefficient of
$\Phi^{-1}$, $j \in {\mathbb N}^0$.

\begin{lemma} {\rm (cf.~\cite[Cor.~3.7]{MiP} for $q=1$)} \label{Lemma5.1}
\newline
   Let $W \in \tilde{{\mathcal W}}_q([-\pi,\pi))$. Then
   \begin{equation} \label{5.2}
      P_{S_2} = I - \sum_{j=0}^n (B_jB_j^*)^{-1} \cdot \sum_{j=0}^n B_j e_j
      (\Phi^*)^{-1}
   \end{equation}
   and
   \begin{equation}   \label{5.3}
      \Delta_{S_2} = \left( \sum_{j=0}^n B_jB_j^* \right)^{-1} \, .
   \end{equation}
\end{lemma}

P r o o f :
By Lemma \ref{Lemma3.8} and Theorems \ref{Theorem3.9} and \ref{Theorem3.11}
the set $S_2$ is $\mathbb Z$-exact. Taking into account Theorem \ref{Theorem3.2}
and the fact that the outer function $\Phi^{-1}$ belongs to $H^2$ as soon as
$W \in \tilde{{\mathcal W}}_q([-\pi,\pi))$, the proof consists of straightforward
calculations which will be omitted.
\hfill $\Box$

\medskip
Our goal is to establish (\ref{5.2}) and (\ref{5.3}) under the weaker assumption
that ${\rm log} \, {\rm det} \, W$ is integrable. This can be done by an approximation
procedure. Our approach is similar to that one presented in the proof of
\cite[Th.~3]{CMiP} for the case $q=1$. However, since we wish to compute not
only the prediction error as done there but also the orthogonal projection, we give
a complete proof of our generalized result.

Let $W_m = W + 1/m \cdot I$, $m \in \mathbb N$. Denote by $(.,.)_m$ the
${\mathbf M}_q$-valued inner product of $L^2(W_m)$, by ${\mathcal M}_m(S_2)$ the
closure of ${\mathcal T}(S_2)$ in $L^2(W_m)$, and by $P_{S_2}^{(m)}$ the
orthogonal projection of $L^2(W_m)$ onto ${\mathcal M}_m(S_2)$. Note that $L^2(W_m)$
can be considered as a subset of $L^2(W)$, and that the inclusion
\begin{equation} \label{5.4}
    {\mathcal M}_m (S_2) \subseteq {\mathcal M}(S_2) \, , \qquad m \in
    \mathbb N \, ,
\end{equation}
holds. Furthermore,
\begin{equation} \label{5.5}
    \lim_{m \to \infty} (I-P_{S_2}^{(m)}I,I-P_{S_2}^{(m)}I)_m =
    (I-P_{S_2}I, I-P_{S_2}I ) \, .
\end{equation}
Let $\Phi_m$ be the corresponding outer factor of $W_m$, i.e.
\[
   W_m=\Phi_m^* \Phi_m \qquad {\rm and} \qquad \int \Phi_m \, d \lambda
   \in {\mathbf M}_q^> \, , \,\, m \in \mathbb N \, .
\]
The following result is an easy consequence of \cite{DLE}, so the proof is omitted.

\begin{lemma} \label{Lemma5.2}
   Let $W \in {\mathcal W}_q'([-\pi,\pi))$. Then $\lim_{m \to \infty} \Phi_m = \Phi$
   w.r.t.~the topology of $L^2$.
\end{lemma}

Let $B_{jm}$ be the $j$-th Taylor coefficient of $\Phi_m^{-1}$. Set
\[
   T_m = \left( \sum_{j=0}^n B_{jm} B_{jm}^* \right)^{-1} \cdot \sum_{j=0}^n
   B_{jm} e_j
\]
and
\[
   T= \left( \sum_{j=0}^n B_j B_j^* \right)^{-1} \cdot \sum_{j=0}^n B_j e_j \, .
\]

\begin{lemma}   \label{Lemma5.3}
   Let $W \in {\mathcal W}_q'([-\pi,\pi))$. Then there exists a subsequence
   $\{ m_l \}_{l \in \mathbb N}$ of $\mathbb N$ such that
   \[
      \lim_{l \to \infty} \left\| I -T (\Phi^*)^{-1} -P_{S_2}^{(m_l)}I
      \right\| = 0 \, .
   \]
\end{lemma}

P r o o f :
Lemma \ref{Lemma5.2} implies that there exists a subsequence $\{ m_l \}_{l
\in \mathbb N}$ of $\mathbb N$ such that $\lim_{l \to \infty} \Phi_{m_l} =
\Phi$ a.e.. Since
\[
   \| (\Phi_{m_l}^*)^{-1} - (\Phi^*)^{-1} \|^2 \leq
   \int | \Phi- \Phi_{m_l} |^2 | \Phi_{m_l}^{-1} |^2 | \Phi^{-1} |^2
   |W^{1/2}|^2 \, d \lambda =
   \int | \Phi-\Phi_{m_l}|^2 |W_{m_l}^{-1} | \, d \lambda
\]
and $|\Phi - \Phi_{m_l} |^2 |W_{m_l}^{-1}| \leq 2 \cdot (|\Phi|^2 + |\Phi_{m_l}|^2)
|W_{m_l}^{-1} | \leq 4$ Lebesgue's dominated convergence theorem yields
\begin{equation} \label{5.6}
    \lim_{l \to \infty} \| (\Phi_{m_l}^*)^{-1} -(\Phi^*)^{-1} \| = 0 \, .
\end{equation}
Applying Lemma \ref{Lemma5.1} we can write
\begin{eqnarray} \label{5.7}
    \| I - T (\Phi   ^*)^{-1} - P_{S_2}^{(m_l)}I \| & = &
       \| T_{m_l}(\Phi_{m_l}^*)^{-1}-T (\Phi^*)^{-1} \| \\ \nonumber
       & \leq &
       \| T_{m_l}(\Phi_{m_l}^*)^{-1} - T_{m_l}(\Phi^*)^{-1} \| +
       \| T_{m_l}(\Phi^*)^{-1} - T (\Phi^*)^{-1} \| \, .  
\end{eqnarray}
Lemma \ref{Lemma5.2} implies that for every $j \in {\mathbb N}^0$ the $j$-th
Taylor coefficient of $\Phi_m$ tends to the $j$-th Taylor coefficient of
$\Phi$ as $m$ tends to infinity. Therefore, $\lim_{l \to \infty} T_{m_l} = T$
uniformly on $[-\pi,\pi)$, and the second summand on the right hand side of
(\ref{5.6}) tends to zero as $l$ tends to infinity. Also the first summand
there tends to zero because of (\ref{5.6}) and the uniform boundedness of the set
$\{ T_{m_l} : l \in \mathbb N \}$.
\hfill $\Box$

\begin{theorem} \label{Theorem5.4}
    Let ${\mathbf G}={\mathbb Z}$, $S_1$ be the set described
    at (\ref{2.2}) and $W \in {\mathcal W}_q'([-\pi,\pi))$. Then
    \begin{equation}  \label{5.8}
       P_{S_2}I = I - \left( \sum_{j=0}^n B_j B_j^* \right)^{-1} \cdot
                  \sum_{j=0}^n B_j e_j (\Phi^*)^{-1}  \, ,
    \end{equation}
    \begin{equation}   \label{5.9}
       \Delta_{S_2} = \left( \sum_{j=0}^n B_j B_j^* \right)^{-1}
    \end{equation}
    and
    \begin{equation} \label{5.10}
       \delta_{S_2} = \left[ {\rm det} \, \left( \sum_{j=0}^n B_jB_j^* \right)
       \right]^{-1/q}     \, .
    \end{equation}
\end{theorem}

P r o o f :
From (\ref{5.5}) and Lemma \ref{Lemma5.3} we get
\[
 \| I-P_{S_2}I \| \geq \lim_{l \to \infty} \| I-P_{S_2}^{(m_l)} I\| =
 \| I -(I-T(\Phi^*)^{-1}) \| \, .
\]
If we take
(\ref{5.4}) into account, we can conclude that $(I-T(\Phi^*)^{-1}) \in
{\mathcal M}(S_2)$. This yields (\ref{5.8}). To obtain (\ref{5.9}) combine
formula (\ref{5.8}) with (\ref{2.1}). Finally, (\ref{5.10}) is a consequence
of (\ref{3.6}).
\hfill $\Box$

\medskip
We conclude our paper by stating the prediction error matrix $\Delta_S$ for
sets $S$ of the form
\begin{equation} \label{2.3}
   S_3 = {\mathbf N}^c \cup \{ n \}
\end{equation}
and
\begin{equation} \label{2.4}
   S_4 = {\mathbf N}^c \setminus \{ -n \}
\end{equation}
for some $n \in \mathbf N$. The univariate versions of the results below can
be found at \cite[Th.~5 and 6]{CMiP}. Let $A_j$
be the $j$-th Taylor coefficient of the outer function $\Phi$ which has the
properties (\ref{5.1}), $j \in {\mathbb N}^0$.

\begin{theorem}  \label{Theorem5.5}
    Let ${\mathbf G}= \mathbb Z$, $S_3$ be the set described at (\ref{2.4}) and
    $W \in {\mathcal W}_q'([-\pi,\pi ))$. Then
    \begin{equation}   \label{5.11}
       \Delta_{S_3} = A_0^* \left( I- A_n \left( \sum_{j=0}^n A_j^*A_j \right)^{-1}
       A_n^* \right) A_0 \, .
    \end{equation}
    If $A_n$ is regular then
    \begin{equation} \label{5.12}
        \Delta_{S_3} = A_0^*(A_n^*)^{-1} \left( \sum_{j=0}^{n-1} A_j^*A_j \right)
                       \left( \sum_{j=0}^{n} A_j^*A_j \right)^{-1} A_n^*A_0 \, .
    \end{equation}
\end{theorem}

We omit the proof since (\ref{5.11}) can be obtained by a straightforward
generalization of the proof given for \cite[Th.~5]{CMiP} to the multivariate
case, and (\ref{5.12}) follows by some simple matrix computations from
(\ref{5.11}).

\begin{theorem} \label{Theorem5.6}
    Let ${\mathbf G}= \mathbb Z$, $S_4$ be the set described at (\ref{2.3}) and
    $W \in {\mathcal W}_q'([-\pi,\pi ))$. Then
    \begin{equation}   \label{5.13}
       \Delta_{S_4} = A_0 \left( I - B_n \left( \sum_{j=0}^n B_j^*B_j \right)^{-1}
       B_n^* \right)^{-1} A_0^* \, .
    \end{equation}
    If $B_n$ is regular then
    \begin{equation} \label{5.14}
       \Delta_{S_4} = A_0 (B_n^*)^{-1} \left( \sum_{j=0}^n B_j^*B_j \right)
                      \left( \sum_{j=0}^{n-1} B_j^*B_j \right)^{-1} B_n^* A_0^* \, .
    \end{equation}
\end{theorem}

P r o o f :
Since the set $S_4$ is $\mathbb Z$-exact by Theorem \ref{Theorem3.11}, Theorem
\ref{Theorem3.9} and Corollary \ref{Corollary3.10} and since the set $S_4^c$ is
the reflection of $S_3$ about the origin, the result immediately follows from
the Theorems \ref{Theorem3.2} and \ref{Theorem5.5} whenever $W \in
\tilde{{\mathcal W}}_q([-\pi,\pi ))$. If merely $W \in {\mathcal W}_q'([-\pi,\pi ))$
approximate $W$ by the sequence $\{ W+ 1/m \cdot I \}_{m \in \mathbb N}$.
\hfill $\Box$

%%%%%%%%%%%%%%%%%%%%%%%%%%%
\renewcommand{\refname}{{\bf References}}

%%%%%%%%%%%%%%%%%%%%%%%%%%%

\bigskip
Universit\"at Leipzig

Fakult\"at f\"ur Mathematik und Informatik

Mathematisches Institut

Augustusplatz 10

D-04109 Leipzig

Fed.~Rep.~Germany

\end{document}